\documentclass[12pt]{article}
\usepackage{cite,amsfonts,amssymb}
\newtheorem{prop}{Proposition}

\newtheorem{lem}{Lemma}

\newtheorem{theor}{Theorem}

\newtheorem{cor}{Corollary}

\newtheorem{rem}{Remark}

\newtheorem{defin}{Definition}

\newcommand{\bgproof}{\noindent {\bf Proof} \hspace{2mm}}
\newcommand{\edproof}{\hfill $\blacksquare$ \vspace{3mm}}

\title{Strong solutions to stochastic Volterra equations}

\author{{\large\sf Anna Karczewska}\\[2mm]
  \normalsize\it
 Department of Mathematics,
 University of Zielona G\'ora\\ \normalsize\it
 ul. Szafrana 4a, 65-246 Zielona G\'ora, Poland,~
 e-mail: A.Karczewska@im.uz.zgora.pl\\[2mm]
 {\large\sf Carlos Lizama}\\[2mm] \normalsize\it
 Universidad de Santiago de Chile, Departamento de
Matem\'atica, Facultad de Ciencias,\\ \normalsize\it Casilla
307-Correo 2,Santiago, Chile,~ e-mail: clizama@lauca.usach.cl }

\begin{document}

\maketitle

\def\thefootnote{}
\footnotetext{\noindent {\em 2000 Mathematics Subject
Classification:}
primary: 60H20; secondary: 60H05, 45D05.\\
{\em Key words and phrases:} stochastic linear Volterra equation,
strong solution resolvent, mild solution, stochastic convolution.}

\begin{abstract}
In this paper stochastic Volterra equations admitting
exponentially bounded resolvents are studied. After obtaining
convergence of resolvents, some properties for stochastic
convolutions are studied. Our main result provide sufficient
conditions for strong solutions to stochastic Volterra equations.

\end{abstract}

\section{Introduction}\label{sSW1}

We deal with the following stochastic Volterra equation in a
separable Hilbert space $H$
\begin{equation} \label{eSW1}
X(t) = X_0 + \int_0^t a(t-\tau)\,AX(\tau)d\tau + \int_0^t
\Psi(\tau)\,dW(\tau)\;, \quad\quad t\geq 0\;,
\end{equation}
where $X_0\in H$,
$a\in L^1_{\mathrm{loc}}(\mathbb{R}_+)$ and $A$ is a closed
unbounded linear operator in $H$ with a dense domain $D(A)$. The
domain $D(A)$ is equipped with the graph norm $|\cdot |_{D(A)}$ of
$A$, i.e.\ $|h|_{D(A)}:=(|h|_H^2+|Ah|_H^2)^{1/2}$, where
$|\cdot|_H$ denotes the norm in $H$.

In this work the equation (\ref{eSW1}) is driven by a cylindrical
Wiener process $W$ and $\Psi$ is an appropriate process defined
later. Let us emphasize that the results obtained for cylindrical
Wiener process $W$ are valid for classical (genuine) Wiener
process, too.

Equation (\ref{eSW1}) arises, in the deterministic case, in a
variety of applications as model problems. Well-known techniques
like localization, perturbation, and coordinate transformation
allow to transfer results for such problems to parabolic
integro-differential equations on smooth domains, see
\cite[Chapter I, Section 5]{Pr2}. In these applications, the
operator $A$ typically is a differential operator acting in
spatial variables, like the Laplacian, the Stokes operator, or the
elasticity operator. The function $a$ should be thought of as a
kernel like $ a(t) = e^{-\eta t} t^{\beta -1}/\Gamma (\beta); \eta
\geq 0, \beta \in (0,2).$ Equation (\ref{eSW1}) is an abstract
stochastic version of the mentioned deterministic model problems.
The stochastic approach to integral equations has been recently
used due to the fact that in applications the level of accuracy
for a given deterministic model not always seem to be
significantly changed with increasing model complexity. Instead,
the stochastic approach provides better results. A typical example
is the use of stochastic integral equations in rainfall-runoff
models, see \cite{Hr-Wh89}.

Our main results concerning (\ref{eSW1}), rely essentially on
techniques using a strongly continuous family of operators $S(t),
t\geq 0$, defined on the Hilbert space $H$ and called the {\tt
resolvent} (precise definition will be given below). Hence, in
what follows, we assume that the deterministic version of equation
(\ref{eSW1}) is {\tt well-posed}, that is, admits a resolvent
$S(t), t\geq 0$. Our aim is to provide sufficient conditions to
obtain a  strong solution to the stochastic equation (\ref{eSW1}).

This paper is organized as follows. In section 2 we prove the main
deterministic ingredient for our construction; this is an
extension of results of \cite{ClNo} allowing here that the
operator $A$ in (\ref{eSW1}) will be the generator of a
$C_0$-semigroup, not necessarily of contraction type. Section 3
contains the main definitions and concepts used in the paper. In
Section 4 we compare mild and weak solutions while in the last
section we provide sufficient condition for stochastic convolution
to be a strong solution to the equation (\ref{eSW1}). We note that
this is an improvement of the known results about existence of
strong solutions for stochastic differential equations.

\section{Convergence of resolvents}\label{sSW2}

In this section we recall some definitions connected with the
deterministic version of the equation (\ref{eSW1}), that is, the
equation
\begin{equation} \label{eSW1d}
u(t) = \int_0^t a(t-\tau)\,Au(\tau)d\tau +f(t), \quad t\geq 0,
\end{equation}
in a Banach space $B$. In (\ref{eSW1d}), the operator $A$ and the
kernel function $a$ are the same as previously considered in the
introduction and $f$ is a $B$-valued function.

Problems of this type have attracted much interest during the last
decades, due to their various applications in mathematical physics
like viscoelasticity, thermodynamics, or electrodynamics with
memory, cf. \cite{Pr2}.

By $S(t),~t\geq 0$, we denote the family of resolvent operators
corresponding to the Volterra equation (\ref{eSW1d}), if it
exists, and defined as follows.

\begin{defin}\label{dSW1} (see, e.g.\ \cite{Pr2})\\
A family $(S(t))_{t\geq 0}$ of bounded linear operators in $B$ is
called {\tt resolvent} for (\ref{eSW1d}) if the following
conditions are satisfied:
\begin{enumerate}
\item $S(t)$ is strongly continuous on $\mathbb{R}_+$ and $S(0)=I$;
\item $S(t)$ commutes with the operator $A$, that is, $S(t)(D(A))\subset
D(A)$ and $AS(t)x=S(t)Ax$ for all $x\in D(A)$ and $t\geq 0$;
\item the following {\tt resolvent equation} holds
\begin{equation} \label{eSW2}
S(t)x = x + \int_0^t a(t-\tau) AS(\tau)x d\tau
\end{equation}
for all $x\in D(A),~t\geq 0$.
\end{enumerate}
\end{defin}

Necessary and sufficient conditions for the existence of the
resolvent family have been studied in \cite{Pr2}. Let us emphasize
that the resolvent $S(t),~t\geq 0$, is determined by the operator
$A$ and the function $a$, so we also say that the pair $(A,a)$
admits a resolvent family. Moreover, as a consequence of the
strong continuity of $S(t)$ we have~ $\sup_{t\leq T} \; ||S(t)|| <
+\infty$ ~for any~ $T\ge 0$.

We shall use the abbreviation $ (a\star f)(t) = \int_0^t
a(t-s)f(s)ds, ~~ t \in [0,T]$, for the convolution of two
functions.

\begin{defin} \label{dSW2}
We say that function $a\in L^1(0,T)$ is {\tt completely positive}
on $[0,T]$ if for any $\mu\geq 0$, the solutions of the
convolution equations
\begin{equation}
s(t) + \mu (a \star s)(t) =1 \quad \mbox{and} \quad r(t) + \mu (a
\star r)(t) = a(t)
\end{equation}
satisfy $s(t)\geq 0$ and $ r(t) \geq 0 $ on $[0,T]$.
\end{defin}

\noindent Kernels with this property have been introduced by
Cl\'{e}ment and Nohel \cite{ClNo}. We note that the class of
completely positive kernels appears naturally in the theory of
viscoelasticity. Several properties and examples of such kernels
appears in \cite[Section 4.2]{Pr2}.

\begin{defin}\label{dSW3}
Suppose $S(t),~t\geq 0$, is a resolvent for (\ref{eSW1d}). $S(t)$
is called {\tt expo\-nen\-tially bounded} if there are constants
$M\geq 1$ and $\omega\in\mathbb{R}$ such that
$$ ||S(t)|| \leq M\,e^{\omega t}, \mbox{~~for all~~} t\geq 0. $$
$(M,\omega)$ is called a {\tt type} of $S(t)$.
\end{defin}

Let us note that in contrary to the case of semigroups, not every
resolvent needs to be exponentially bounded even if the kernel
function $a$ belongs to $L^1(\mathbb{R}_+)$. The resolvent
 version of the Hille--Yosida theorem (see, e.g.,
\cite[Theorem 1.3]{Pr2}) provides the class of equations that
admit exponentially bounded resolvents. An important class of
kernels providing such class of resolvents are $a(t)=
t^{\beta-1}/\Gamma(\beta),~\alpha\in (0,2)$ or the class of
completely monotonic functions. For details, counterexamples and
comments we refer to \cite{DePr}.

In this paper the following result concerning convergence of
resolvents for the equation (\ref{eSW1}) in a Banach space $B$
will play the key role. It corresponds to a generalization of the
results of Cl\'{e}ment and Nohel obtained in \cite{ClNo} for
contraction semigroups.

\begin{theor} \label{pSW2}
Let $A$ be the generator of a $C_0$-semigroup in $B$ and suppose
the kernel function $a$ is completely positive. Then $(A,a)$
admits an exponentially bounded resolvent $S(t)$. Moreover, there
exist bounded operators $A_n$ such that $(A_n,a)$ admit resolvent
families $S_n(t)$ satisfying $ ||S_n(t) || \leq Me^{w_0 t}~ (M\geq
1,~w_0\geq 0)$ for all $t\geq 0$ and
\begin{equation} \label{eSW4}
S_n(t)x \to S(t)x \quad \mbox{as} \quad n\to +\infty
\end{equation}
for all $x \in B,\; t\geq 0.$ Additionally, the convergence is
uniform in $t$ on every compact subset of $ \mathbb{R}_+$.
\end{theor}

\bgproof The first assertion follows directly from \cite[Theorem
5]{Pr1} (see also \cite[Theorem 4.2]{Pr2}). Since $A$ generates a
$C_0$-semigroup $T(t),~t\geq 0$, the resolvent set $\rho(A) $ of
$A$ contains the ray $ [w,\infty)$ and
$$
||R(\lambda,A)^k || \leq \frac{M}{(\lambda -w)^k } \qquad
\mbox{for } \lambda > w, \qquad k\in \mathbb{N}.
$$

Define
\begin{equation} \label{eSW5}
A_n := n AR(n,A) = n^2 R(n,A) - nI, \qquad n> w
\end{equation}
the {\it Yosida approximation} of $A$.

Then
\begin{eqnarray*}
||e^{t A_n} || &=& e^{-nt} || e^{n^2 R(n,A)t} || \leq
e^{-nt} \sum_{k=0}^{\infty} \frac{n^{2k} t^k}{k!} ||R(n,A)^k|| \\
&\leq& M e^{(-n + \frac{n^2}{n-w})t} = M e^{ \frac{nwt}{ n-w}}.
\end{eqnarray*}
Hence, for $n > 2w$ we obtain
\begin{equation}{\label{eSW6}}
|| e^{A_n t} || \leq M e^{2wt}.
\end{equation}
Taking into account the above estimate and the complete positivity
of the kernel function $a$, we can follow the same steps as in
\cite[Theorem 5]{Pr1} to obtain that there exist constants $M_1
>0 $ and $ w_1 \in \mathbb{R} $ (independent of $n$, due to
(\ref{eSW6})) such that
$$
||[H_n(\lambda)]^{(k)} || \leq \frac{M_1}{ (\lambda - w_1)^{k+1}}
\quad \mbox{ for } \lambda > w_1,
$$
where $H_n(\lambda):= (\lambda - \lambda \hat
a(\lambda)A_n)^{-1}.$ Here and in the sequel the hat indicates the
Laplace transform. Hence, the generation theorem for resolvent
families implies that for each $ n > 2\omega$, the pair $(A_n,a)$
admits resolvent family $S_n(t)$ such that
\begin{equation}{\label{eSW7}}
||S_n(t) || \leq M_1 e^{w_1 t}.
\end{equation}
In particular, the Laplace transform $ \hat S_n(\lambda) $ exists
and satisfies
$$
\hat S_n(\lambda ) = H_n(\lambda) = \int_0^{\infty} e^{-\lambda t}
S_n(t)dt, \qquad \lambda > w_1.
$$
Now recall from semigroup theory that for all $\mu$ sufficiently
large we have
$$ R(\mu,A_n)= \int_0^\infty e^{-\mu t} \,e^{A_nt}\,dt $$
as well as,
$$ R(\mu,A)= \int_0^\infty e^{-\mu t} \,T(t)\,dt\,. $$

Since $\hat a(\lambda) \to 0$ as $ \lambda \to \infty$, we deduce
that for all $\lambda$ sufficiently large, we have
$$ 
H_n(\lambda) := \frac{1}{\lambda \hat a (\lambda) } R(
\frac{1}{\hat a (\lambda) }, A_n) = \frac{1}{\lambda \hat a
(\lambda)} \int_0^{\infty} e^{(-1/\hat a(\lambda))t} e^{A_n t}
dt\,,
$$ 
and
$$ 
H(\lambda) := \frac{1}{\lambda \hat a (\lambda) } R( \frac{1}{\hat
a (\lambda) }, A) = \frac{1}{\lambda \hat a (\lambda)}
\int_0^{\infty} e^{(-1/\hat a(\lambda))t} T(t) dt\,.
$$ 

Hence, from the identity
$$ H_n(\lambda) - H(\lambda) = \frac{1}{\lambda \hat a (\lambda) }
[R(\frac{1}{\hat a (\lambda) }, A_n)-R(\frac{1}{\hat a (\lambda)
}, A)]
$$
and the fact that $R(\mu,A_n)\to R(\mu,A)$ as $n\to\infty$ for all
$\mu$ sufficiently large (see, e.g.\ \cite[Lemma~7.3]{Pa}, we
obtain that
\begin{equation}{ \label{eSW8}}
H_n(\lambda) \to H(\lambda) \quad \mbox{as } n \to \infty\;.
\end{equation}
Finally, due to (\ref{eSW7}) and (\ref{eSW8}) we can use the
Trotter-Kato theorem for resolvent families of operators (cf.
\cite[Theorem 2.1]{Li}) and the conclusion follows. \edproof

An analogous result like Theorem \ref{pSW2} holds in other cases.

\begin{theor} \label{pSW2a}
 Let $A$ be the generator of a strongly continuous cosine family.
 Suppose any of the following:
 \begin{description}
  \item[~~(i)] $a\in L_\mathrm{loc}^1(\mathbb{R}_+)$ is completely
  positive;
  \item[~(ii)] the kernel fuction $a$ is a creep function with $a_1$
  log-convex;
  \item[(iii)] $a=c\star c$ with some completely positive
  $c\in L_\mathrm{loc}^1(\mathbb{R}_+)$.
 \end{description}
Then $(A,a)$ admits an exponentially bounded resolvent $S(t)$.
Moreover, there exist bounded operators $A_n$ such that $(A_n,a)$
admit resolvent families $S_n(t)$ satisfying $ ||S_n(t) || \leq
Me^{w_0 t}$ $(M\geq 1,~w_0\geq 0)$ for all $t\geq 0,~n\in
\mathbb{N},$ and
$$ 
S_n(t)x \to S(t)x \quad \mbox{as} \quad n\to +\infty
$$ 
for all $x \in B,\; t\geq 0.$ Additionally, the convergence is
uniform in $t$ on every compact subset of $ \mathbb{R}_+$.
\end{theor}

The proof follows from \cite[Theorem 4.3]{Pr2}, where the
definition of a creep function can be found, or \cite[Theorem
6]{Pr1} and proof of Theorem \ref{pSW2}. Therefore it is omitted.

\begin{rem}{\label{'comment'}}
 {\em Other examples of the convergence (\ref{eSW4}) for
the resolvents are given, e.g., in \cite{ClNo} and \cite{Fr}. In
the first paper, the operator $A$ generates a linear continuous
contraction semigroup. In the second of the mentioned papers, $A$
belongs to some subclass of sectorial operators and the kernel $a$
is an absolutely continuous function fulfilling some technical
assumptions.}
\end{rem}

\begin{prop}\label{AScom}
 Let $A, A_n$ and $ S_n( t)$ be given as in Theorem \ref{pSW2}. Then $ S_n(t)$ commutes with  the operator $A$, for every
 $n$ sufficiently large and $t\ge 0$.
\end{prop}

\bgproof For each $n$ sufficiently large the bounded operators
$A_n$ admit a resolvent family $S_n(t)$, so by the complex
inversion formula for the Laplace transform we have
$$S_n(t)=\frac{1}{2\pi i}\int_{\Gamma_n} e^{\lambda t}H_n(\lambda) d\lambda$$
where $\Gamma_n$ is a simple closed rectifiable curve surrounding
the spectrum of $A_n$ in the positive sense.

 On the other hand, $H_n(\lambda) := (\lambda - \lambda\hat
a (\lambda) A_n)$ where $ A_n := nA(n -A)^{-1},$ so each $A_n$
commutes with $A$ on $D(A)$ and then each $H_n(\lambda)$ commutes
with $A$, on $D(A)$, too.

Finally, because $A$ is closed, and all the following integrals
are convergent (exist) we have for all $n$ sufficiently large and
$x \in D(A)$
\begin{eqnarray*}
AS_n(t)x & = & A \int_{\Gamma_n} e^{\lambda t} H_n(\lambda)x
d\lambda
 = \int_{\Gamma_n} e^{\lambda t} AH_n(\lambda)x d\lambda \\
 & = & \int_{\Gamma_n} e^{\lambda t} H_n(\lambda)Ax d\lambda
  = S_n(t) Ax \,.
\end{eqnarray*}

\edproof

\section{Solutions and the stochastic convolution} \label{sSW3}

Let $H$ and $U$ be two separable Hilbert spaces and $Q\in L(U)$ be
a linear bounded symmetric nonnegative operator. A Wiener process
$W$ with covariance operator $Q$ is defined on a probability space
$(\Omega,\mathcal{F},(\mathcal{F}_t)_{t\geq 0},P)$. We assume that
the process $W$ is a cylindrical one, that is, we do not assume
that $\mathrm{Tr} Q<+\infty$. In this case, the process $W$ has
values in some superspace of $U$. Let us note that the results
obtained in the paper for cylindrical Wiener process are valid in
classical Wiener process, too. Namely, when $\mathrm{Tr}
Q<+\infty$, we can take $U=H$ and $\Psi=I$.

This is apparently well-known that the construction of the
stochastic integral with respect to cylindrical Wiener process
requires some particular terms. We will need the subspace
$U_0:=Q^{1/2}(U)$ of the space $U$, which endowed with the inner
product $\langle u,v\rangle_{U_0}:=\langle Q^{-1/2}u,
Q^{-1/2}v\rangle_U$ forms a Hilbert space. Among others, an
important role is played by the space of Hilbert-Schmidt
operators. The set $L_2^0:=L_2(U_0,H)$ of all Hilbert-Schmidt
operators from $U_0$ into $H$, equipped with the norm
$|C|_{L_2(U_0,H)}:=(\sum_{k=1}^{+\infty} |C f_k|_H^2)^{1/2}$,
where $\{f_k\}$ is an orthonormal basis of $U_0$, is a separable
Hilbert space. \label{H-Sn}

According to the theory of stochastic integral with respect to
cylindrical Wiener process we have to assume that $\Psi$ belongs
to the class of measurable $L_2^0$-valued processes.

Let us introduce the norms
\begin{eqnarray*}
||\Psi||_t &:=& \left\{\mathbb{E}\left( \int_0^t
|\Psi(\tau)|_{L_2^0}^2\,d\tau \right) \right\}^{\frac{1}{2}} \\
& = & \left\{\mathbb{E} \int_0^t \left[ \mathrm{Tr}
(\Psi(\tau)Q^{\frac{1}{2}}) (\Psi(\tau)Q^{\frac{1}{2}})^* \right]
d\tau \right\}^{\frac{1}{2}}, \quad t\in [0,T].
\end{eqnarray*}
By $\mathcal{N}^2(0,T;L_2^0)$ we denote a Hilbert space of all
$L_2^0$-predictable processes $\Psi$ such that $|| \Psi||_T <
+\infty$.

It is possible to consider a more general class of integrands,
see, e.g.\ \cite{LiSh}, but in our opinion it is not worthwhile to
study the general case. That case produces a new level of
difficulty additionally to problems related to long time memory of
the system. So, we shall study the equation (\ref{eSW1}) under the
below {\sc Probability Assumptions} (abbr.\ (PA)):
\begin{enumerate}
\item $X_0$ is an $H$-valued, $\mathcal{F}_0$-measurable random variable;
\item $\Psi\in \mathcal{N}^2(0,T;L_2^0)$ and the interval $[0,T]$ is fixed.
\end{enumerate}

\begin{defin} \label{dSW4}
Assume that (PA) hold. An $H$-valued predictable process
$X(t),~t\in [0,T]$, is said to be a ~{\tt strong solution}~ to
(\ref{eSW1}), if $X$ has a version such that $P(X(t)\in D(A))=1$
for almost all $t\in [0,T]$; for any $t\in [0,T]$
\begin{equation} \label{eSW3.1}
\int_0^t |a(t-\tau)AX(\tau)|_H \,d\tau<+\infty,\quad P-a.s.
\end{equation}
and for any $t\in [0,T]$ the equation (\ref{eSW1}) holds $P$-a.s.
\end{defin}

Let $A^*$ denote the adjoint of $A$ with a dense domain
$D(A^*)\subset H$ and the graph norm $|\cdot |_{D(A^*)}$.

\begin{defin} \label{dSW5}
Let (PA) hold. An $H$-valued predictable process $X(t),~t\in
[0,T]$, is said to be a {\tt weak solution} to (\ref{eSW1}), if
$P(\int_0^t|a(t-\tau)X(\tau)|_H d\tau<+\infty)=1$ and if for all
$\xi\in D(A^*)$ and all $t\in [0,T]$ the following equation holds
$$
\langle X(t),\xi\rangle_H = \langle X_0,\xi\rangle_H + \langle
\int_0^t a(t-\tau)X(\tau)\,d\tau, A^*\xi\rangle_H + \langle
\int_0^t \Psi(\tau)dW(\tau),\xi\rangle_H, ~~P\mathrm{-a.s.}
$$
\end{defin}

\begin{rem}{\em This definition has sense for a cylindrical
Wiener process because the scalar process $\langle \int_0^t
\Psi(\tau)dW(\tau),\xi\rangle_H$, for $t\in [0,T]$, is
well-defined.}
\end{rem}

\begin{defin} \label{dSW6}
Assume that $X_0$ is an $H$-valued $\mathcal{F}_0$-measurable
random variable. An $H$-valued predictable process $X(t),~t\in
[0,T]$, is said to be a {\tt mild solution} to the stochastic
Volterra equation (\ref{eSW1}), if
\begin{equation}\label{edf6}
 \mathbb{E}\left( \int_0^t |S(t-\tau)\Psi(\tau)|_{L_2^0}^2 \,d\tau
 \right) <+\infty
 \quad \mbox{for} \quad t\leq T
\end{equation}
  and, for arbitrary $t\in [0,T]$,
\begin{equation}\label{eSW9}
X(t) = S(t)X_0 + \int_0^t S(t-\tau)\Psi(\tau)\,dW(\tau), \quad
P-a.s.
\end{equation}
where $S(t)$ is the resolvent for the equation (\ref{eSW1d}), if
it exists.
\end{defin}

In some cases weak solutions of equation (\ref{eSW1}) coincides
with mild solutions of (\ref{eSW1}), see e.g.\ \cite{Ka}. In
consequence, having results for the convolution on the right hand
side of (\ref{eSW9}) we obtain results for weak solutions.

In the paper we will use the following well-known result.

\begin{prop} \label{pSW1} (see, e.g.\cite[Proposition 4.15]{DaPrZa})\\
Assume that $A$ is a closed linear unbounded operator with the
dense domain $D(A)\subset H$ and $\Phi(t),~t\in [0,T]$ is an
$L_2(U_0,H)$-predictable process. If $~\Phi(t)(U_0)\subset D(A), ~P-a.s.$
for all $t\in [0,T]$ and
$$ P\left( \int_0^T | \Phi(s)|_{L_2^0}^2\,ds <\infty \right) =1,~~
P\left( \int_0^T |A \Phi(s)|_{L_2^0}^2\,ds <\infty \right) =1,
$$
then 
$$  P\left( \int_0^T \Phi(s)\,dW(s) \in D(A)
\right) =1 \mbox{~~and~~}
 A \int_0^T \Phi(s)\,dW(s) = \int_0^T A\Phi(s)\,dW(s), ~~P-a.s.
$$
\end{prop}


In what follows we assume that (\ref{eSW1d}) admits a resolvent
family $S(t)$, $t \geq 0.$ We introduce the stochastic convolution
\begin{equation} \label{eSW18a}
W^\Psi(t) := \int_0^t S(t-\tau)\Psi(\tau)\,dW(\tau),
\end{equation}
where $\Psi$ belongs to the space $\mathcal{N}^2(0,T;L_2^0)$. Note
that, because resolvent operators $S(t),~t\geq 0$, are bounded,
then $S(t-\cdot)\Psi (\cdot)\in \mathcal{N}^2(0,T;L_2^0)$, too.

Let us formulate some auxiliary results concerning the convolution
$W^\Psi (t)$.

\begin{prop}\label{pr2proper}
Assume that (\ref{eSW1d}) admits resolvent operators $S(t),~t\geq
0$. Then, for arbitrary process $\Psi\in\mathcal{N}^2(0,T;L_2^0)$,
the process $W^\Psi(t),~t\geq 0$, given by (\ref{eSW18a}) has a
predictable version.
\end{prop}

\begin{prop}\label{pr3proper}
Assume that $\Psi\in\mathcal{N}^2(0,T;L_2^0)$. Then the process
$W^\Psi(t),~t\geq 0$, defined by (\ref{eSW18a}) has square
integrable trajectories.
\end{prop}

For  the proofs of Propositions \ref{pr2proper} and
\ref{pr3proper} we refer to \cite{Ka}.

\begin{prop} \label{pr4proper}
Let $a\in BV(\mathbb{R}_+)$ and suppose that (\ref{eSW1d}) admits
a resolvent family $S \in C^1(0,\infty; L(H)).$ Let $X$ be a
predictable process with integrable trajectories. Assume that $X$
has a version such that $P(X(t)\in D(A))=1$ for almost all $t\in
[0,T]$ and (\ref{edf6}) holds. If for any $t\in [0,T]$ and $\xi\in
D(A^*)$
\begin{equation}\label{deq9}
 \langle X(t),\xi\rangle_H = \langle X_0,\xi\rangle_H +
 \int_0^t \langle
 a(t-\tau)X(\tau),A^*\xi\rangle_H d\tau
 + \int_0^t \langle\xi,\Psi(\tau) dW(\tau)\rangle_H, ~~P-a.s.,
\end{equation}
then
\begin{equation}\label{deq9a}
X(t) = S(t)X_0 +
 \int_0^t S(t-\tau) \Psi (\tau) dW(\tau), \quad t\in[0,T].
\end{equation}
\end{prop}
\bgproof For simplicity we omit the index $_H$ in the inner
product. Since $a\in BV(\mathbb{R}_+),$  we can see that
(\ref{deq9}) implies
\begin{eqnarray}\label{deq14}
 \langle X(t),\xi(t)\rangle &=& \langle X_0,\xi(0)\rangle +
 \int_0^t \langle (\dot{a}\star X)(\tau) +
 a(0)X(\tau),A^*\xi(\tau)\rangle d\tau \nonumber \\
 &+& \int_0^t \langle \Psi(\tau) dW(\tau),\xi(\tau)\rangle
 + \int_0^t \langle X(\tau),\dot{\xi}(\tau) \rangle d\tau, \quad
 \mathrm{~P-a.s.}
\end{eqnarray}
for any $\xi\in C^1([0,t],D(A^*))$ and $t\in [0,T]$. (For details,
see \cite{Ka}).

Now, let us take $\xi(\tau):=S^*(t-\tau)\zeta$ with $\zeta\in
D(A^*)$,
 $\tau\in [0,t]$.
The equation (\ref{deq14}) may be written like
\begin{eqnarray*}
 \langle X(t),S^*(0)\zeta\rangle & = & \langle X_0,S^*(t)\zeta\rangle
 + \int_0^t \langle (\dot{a}\star X)(\tau)
 + a(0)X(\tau),A^*S^*(t-\tau)\zeta\rangle d\tau \nonumber \\
  &+& \int_0^t \langle \Psi(\tau) dW(\tau),S^*(t-\tau)\zeta\rangle
 + \int_0^t \langle X(\tau),(S^*(t-\tau)\zeta)'\rangle d\tau,
\end{eqnarray*}
where the derivative ()' in the last term is taken over $\tau$.

Next, using $S^*(0)=I$, we rewrite
\begin{eqnarray} \label{deq15}
 \langle X(t),\zeta\rangle &=& \langle S(t)X_0,\zeta\rangle +
 \int_0^t \langle S(t-\tau)A\left[\int_0^\tau
 \dot{a}(\tau-\sigma)X(\sigma)d\sigma +a(0)X(\tau)\right],\zeta\rangle
 d\tau \nonumber\\
   &+& \int_0^t \langle S(t-\tau)\Psi(\tau)dW(\tau),\zeta\rangle +
   \int_0^t \langle \dot{S}(t-\tau)X(\tau),\zeta\rangle d\tau.
\end{eqnarray}
To prove (\ref{deq9a}) it is enough to show that the sum of the
first integral and the third one in the equation (\ref{deq15})
gives zero.

Since $S \in C^1(0,\infty; L(H)),$ we can use properties of
resolvent operators and the derivative $\dot{S}(t-\tau)$ with
respect to~$\tau$. Then
\begin{eqnarray*}
 I &:=& \left\langle \int_0^t \dot{S}(t-\tau)X(\tau) d\tau,\zeta \right\rangle =
  \left\langle -\int_0^t \dot{S}(\tau)X(t-\tau) d\tau ,\zeta \right\rangle \\
 &=& \left\langle -\left(
  \int_0^t\left[ \int_0^\tau \dot{a}(\tau-s)AS(s)ds\right] X(t-\tau) d\tau
 - \int_0^t a(0)AS(\tau) X(t-\tau) d\tau \right) ,\zeta \right\rangle \\[2mm]
 & = & \langle -([A(\dot{a}\star S)(\tau) \star X](t)+a(0)A(S\star X)(t))
 ,\zeta\rangle .
\end{eqnarray*}
Note that $a\in BV(\mathbb{R}_+)$ and hence the convolution
$(a\star S)(\tau)$ has sense (see \cite[Section 1.6]{Pr2}).

Since
$$ \int_0^t \langle a(0)AS(t-\tau)X(\tau),\zeta\rangle d\tau =
  \int_0^t \langle a(0)AS(\tau)X(t-\tau),\zeta\rangle d\tau $$
 and
\begin{eqnarray*}
J &:=& \int_0^t \langle S(t-\tau)A\left[\int_0^\tau
 \dot{a}(\tau-\sigma)X(\sigma)d\sigma\right],\zeta\rangle d\tau =
 \int_0^t \langle AS(t-\tau)(\dot{a}\star X)(\tau),\zeta\rangle d\tau = \\[2.5mm]
 &\!\!=\!\!& \langle A(S\star (\dot{a}\star X)(\tau))(t),\zeta\rangle =
 \langle A((S\star\dot{a})(\tau)\star X)(t),\zeta\rangle
 \quad \mbox{for~any~~} \zeta\in D(A^*),
\end{eqnarray*}
so $J=-I$, hence $J+I=0$. This means that (\ref{deq9a}) holds for
any $\zeta\in D(A^*)$. Since $D(A^*)$ is dense in $H^*$, then
(\ref{deq9a}) holds. \edproof

\begin{rem}
{\em If  (\ref{eSW1}) is parabolic and the kernel $a$ is
3-monotone, understood in the sense defined by Pr\"uss
\cite[Section 3]{Pr2}, then $S \in C^1(0,\infty; L(H))$ and $a\in
BV( \mathbb{R}_+)$ respectively.}
\end{rem}

\begin{prop} \label{pr5prop}
Assume that $A$ is a closed linear unbounded operator with the dense
domain $D(A)$, $a\in L_\mathrm{loc}^1(\mathbb{R}_+)$ and
$S(t)$, $t\ge 0$, are resolvent operators for the equation (\ref{eSW1d}).
If $\Psi\in\mathcal{N}^2(0,T;L_2^0)$, then the stochastic convolution
$W^\Psi$ fulfills the equation (\ref{deq9}) with $X_0\equiv 0$.
\end{prop}
\bgproof Let us notice that the process $W^\Psi$ has integrable
trajectories. For any $\xi\in D(A^*)$ we have from (\ref{eSW18a})
\begin{eqnarray*}
 \int_0^t \langle a(t-\tau)W^\Psi(\tau),A^*\xi\rangle_Hd\tau
 &\equiv& \int_0^t \langle a(t-\tau) \int_0^\tau S(\tau-\sigma)\Psi(\sigma)
 dW(\sigma),A^*\xi\rangle_Hd\tau.
 \end{eqnarray*}
 Hence from~Dirichlet's~formula and the stochastic Fubini's
 theorem we get
 \begin{eqnarray*}
 \int_0^t \langle a(t-\tau)W^\Psi(\tau),A^*\xi\rangle_Hd\tau
&=& \int_0^t \langle\left[\int_\sigma^t
a(t-\tau)S(\tau-\sigma)d\tau\right]
 \Psi(\sigma) dW(\sigma),A^*\xi\rangle_H\\
 &=& \langle \!\int_0^t\! \left[\int_0^{t-\sigma}\! a(t-\sigma-z)S(z)dz\right]
 \!\Psi(\sigma) dW(\sigma),A^*\xi\rangle_H.
\end{eqnarray*}
Taking  $z:=\tau-\sigma$ ~and~from definition~of~convolution we
have
\begin{eqnarray*}
 \int_0^t \langle a(t-\tau)W^\Psi(\tau),A^*\xi\rangle_Hd\tau &=& \langle \int_0^t A [ (a\star S)(t-\sigma)]\Psi(\sigma)
dW(\sigma),
  \xi\rangle_H
  \end{eqnarray*}
From~the~resolvent~equation~(\ref{eSW2}) and because $
     A(a\star S)(t-\sigma)x = (S(t-\sigma)-I)x, $ \\
   ~where $ x\in D(A),$ we obtain
   \begin{eqnarray*}
\int_0^t \langle a(t-\tau)W^\Psi(\tau),A^*\xi\rangle_Hd\tau  &=& \langle \int_0^t [S(t-\sigma)-I]\Psi(\sigma) dW(\sigma),\xi\rangle_H = \\
 = \langle \int_0^t S(t-\sigma)\Psi(\sigma) dW(\sigma),\xi\rangle_H
 &-&\langle \int_0^t \Psi(\sigma) dW(\sigma),\xi\rangle_H .
\end{eqnarray*}
Hence, we obtained the following equation
$$
 \langle W^\Psi(t),\xi\rangle_H = \int_0^t \langle a(t-\tau)W^\Psi(\tau),
 A^*\xi\rangle_Hd\tau + \int_0^t \langle \xi,\Psi(\tau)dW(\tau)\rangle_H
$$
for any $\xi\in D(A^*)$. \edproof

\begin{cor} \label{c1prop}
Assume that $A$ is a linear bounded operator in $H$, 
$a\in L_\mathrm{loc}^1(\mathbb{R}_+)$ and $S(t)$, $t\ge 0$, 
are resolvent operators for the equation (\ref{eSW1d}).
If $\;\Psi$ belongs to $\mathcal{N}^2(0,T;L_2^0)$ then
\begin{equation} \label{deq16}
 W^\Psi(t) = \int_0^t a(t-\tau)AW^\Psi(\tau)d\tau
 + \int_0^t \Psi(\tau) dW(\tau)\,.
\end{equation}
\end{cor}

\begin{rem} {\em The formula (\ref{deq16}) says that the convolution $W^\Psi$
is a strong solution to (\ref{eSW1}) with $X_0\equiv 0$ if the
operator $A$ is bounded.}
\end{rem}

\section{Strong solution}\label{SS5}

In this section we provide sufficient conditions under which the
stochastic convolution $W^\Psi (t)$, $t\ge 0$, defined by
(\ref{eSW18a}) is a strong solution to the equation (\ref{eSW1}).

\begin{lem} \label{pSW5}
Let $A$ be a closed linear unbounded operator the with dense
domain $D(A)$ equipped with the graph norm $|\cdot|_{D(A)}$.
Suppose that assumptions of Theorem \ref{pSW2} or Theorem
\ref{pSW2a} hold. If $\Psi$ and $A\Psi$ belong to
$\mathcal{N}^2(0,T;L_2^0)$ and in addition
$\Psi(\cdot,\cdot)(U_0)\subset D(A),$ P-a.s., then (\ref{deq16})
holds.
\end{lem}
\bgproof Because formula (\ref{deq16}) holds for any bounded
operator, then it holds for the Yosida approximation $A_n$ of the
operator $A$, too, that is
$$ W_n^\Psi(t) =
\int_0^t a(t-\tau) A_n W_n^\Psi(\tau)d\tau +
\int_0^t\Psi(\tau)dW(\tau), $$ where
$$ W_n^\Psi(t) := \int_0^t S_n(t-\tau)\Psi(\tau)dW(\tau)$$
and
$$ A_n W_n^\Psi(t) =
A_n \int_0^t S_n(t-\tau)\Psi(\tau)dW(\tau).
$$
Recall that by assumption $\Psi\in \mathcal{N}^2(0,T;L_2^0)$.
Because the operators $S_n(t)$ are deterministic and bounded for
any $t\in [0,T]$, $n\in\mathbb{N}$, then the operators
$S_n(t-\cdot )\Psi(\cdot)$ belong to $\mathcal{N}^2(0,T;L_2^0)$,
too. In consequence, the difference
\begin{equation}\label{eq21b}
 \Phi_n(t-\cdot ) := S_n(t-\cdot )\Psi(\cdot)
  - S(t-\cdot )\Psi(\cdot)
\end{equation}
belongs to $\mathcal{N}^2(0,T;L_2^0)$ for any $t\in [0,T]$ and
$n\in\mathbb{N}$. This means that
\begin{equation}\label{eq22}
 \mathbb{E}\left(\int_0^t |\Phi_n(t-\tau)|_{L_2^0}^2d\tau \right)
 < +\infty
\end{equation}
for any $t\in [0,T]$.

Let us recall that the cylindrical Wiener process $W(t)$, $t\ge
0$, can be written in the form
\begin{equation}\label{eq23}
 W(t) =\sum_{j=1}^{+\infty} f_j\,\beta_j(t),
\end{equation}
where $\{f_j\}$ is an orthonormal basis of $U_0$ and $\beta_j(t)$
are independent real Wiener processes. From (\ref{eq23}) we have
\begin{equation}\label{eq24}
 \int_0^t \Phi_n(t-\tau)\,dW(\tau) = \sum_{j=1}^{+\infty}
 \int_0^t \Phi_n(t-\tau)\,f_j\,d\beta_j(\tau).
\end{equation}
Then, from (\ref{eq22})
\begin{equation}\label{eq25}
 \mathbb{E}\left[\int_0^t \left( \sum_{j=1}^{+\infty}
 |\Phi_n(t-\tau)\,f_j|_H^2 \right) d\tau \right]
 < +\infty
\end{equation}
for any $t\in [0,T]$. Next, from (\ref{eq24}), properties of
stochastic integral and (\ref{eq25}) we obtain for any
$t\in[0,T]$,
\begin{eqnarray*}
 \mathbb{E}\left| \int_0^t \Phi_n(t-\tau)\,dW(\tau) \right|_H^2
 &=& \mathbb{E}\left| \sum_{j=1}^{+\infty}\int_0^t
 \Phi_n(t-\tau)\,f_j\,d\beta_j(\tau) \right|_H^2 \le \\
  \mathbb{E}\left[ \sum_{j=1}^{+\infty} \int_0^t
 |\Phi_n(t-\tau)\,f_j|_H^2 d\tau \right]
 &\le & \mathbb{E}\left[ \sum_{j=1}^{+\infty} \int_0^T
  |\Phi_n(T-\tau)\,f_j|_H^2 d\tau \right] <+\infty.
\end{eqnarray*}

By Theorem \ref{pSW2} or \ref{pSW2a}, the convergence (\ref{eSW4})
of resolvent families is uniform in $t$ on every compact subset of
$\mathbb{R}_+$, particularly on the interval $[0,T]$. Now, we use
(\ref{eSW4}) in the Hilbert space $H$, so (\ref{eSW4}) holds for
every $x\in H$. Then, for any fixed $j$,
\begin{equation}\label{eq26}
 \int_0^T |[S_n(T-\tau)-S(T-\tau)]\,
 \Psi(\tau)\,f_j|_H^2 d\tau
\end{equation}
tends to zero for $n\to +\infty$. Summing up our considerations,
particularly using (\ref{eq25}) and (\ref{eq26}) we can write
\begin{eqnarray*}
\sup_{t\in [0,T]} \, \mathbb{E}\left| \int_0^t
\Phi_n(t-\tau)dW(\tau)
 \right|_H^2 \equiv \sup_{t\in [0,T]}\, \mathbb{E}\left| \int_0^t
 [S_n(t-\tau)-S(t-\tau)] \Psi(\tau)dW(\tau)\right|_H^2
 & \le & \\ \le
 \mathbb{E}\left[ \sum_{j=1}^{+\infty} \int_0^T
 | [ S_n(T-\tau)-S(T-\tau)]\Psi(\tau)\,f_j |_H^2
 d\tau \right] & \to & 0
\end{eqnarray*}
as $n \to +\infty$.

Hence, by the Lebesgue dominated convergence theorem
\begin{equation}\label{eq27}
 \lim_{n\to +\infty} \sup_{t\in [0,T]} \mathbb{E} \left|
 W_n^\Psi(t)- W^\Psi(t)\right|_H^2 =0.
\end{equation}

By assumption, $\Psi(\cdot,\cdot)(U_0)\subset D(A), ~P-a.s.~$
Because $~S(t)(D(A))\subset D(A)$, then 
$S(t-\tau)\Psi(\tau)(U_0) \subset D(A),
~P-a.s.$, for any $\tau\in [0,t],~t\ge 0$.
Hence, by Proposition \ref{pSW1}, $P(W^\Psi (t)\in D(A))=1$. 

For any $n\in\mathbb{N}$, $t\ge 0$, we have
$$ |A_n W_n^\Psi (t) - A W^\Psi (t)|_H \le
    N_{n,1}(t) + N_{n,2}(t), $$
where
\begin{eqnarray*}
 N_{n,1}(t) & := & |A_n W_n^\Psi (t) - A_n W^\Psi (t)|_H , \\
 N_{n,2}(t) & := & |A_n W^\Psi (t) - A W^\Psi (t)|_H =
             |(A_n-A)W^\Psi (t)|_H \,.
\end{eqnarray*}
Then
\begin{eqnarray}\label{eq28}
 |A_n W_n^\Psi (t) - A W^\Psi (t)|_H^2 & \le &
 N_{n,1}^2 (t) + 2 N_{n,1}(t) N_{n,2}(t) + N_{n,2}^2(t) \nonumber \\
 & < & 3[N_{n,1}^2 (t)+N_{n,2}^2(t)].
\end{eqnarray}

Let us study the term $N_{n,1}(t)$. Note that the unbounded
operator $A$ generates a semigroup. Then we have for the Yosida
approximation the following properties:
\begin{equation}\label{eq30}
 A_nx=J_nAx \quad \mbox{for~any~} x\in D(A), \quad \sup_n ||J_n||
 < \infty
\end{equation}
where $A_nx=nAR(n,A)x=AJ_nx$ for any $x\in H$, with
$J_n:=nR(n,A).$ Moreover (see \cite[Chapter II, Lemma 3.4]{EnNa}):
\begin{eqnarray}
 \lim_{n\to\infty}J_nx &=& x \qquad \mbox{for~any~} x\in H, \nonumber \\
 \lim_{n\to\infty} A_nx &=& Ax \qquad \mbox{for~any~} x\in D(A).
 \label{eq31}
\end{eqnarray}
By Proposition \ref{AScom}, $AS_n(t)x = S_n(t)Ax$ for all $x \in
D(A)$. So, by Propositions \ref{AScom} and \ref{pSW1} and the 
closedness of $A$ we can write
\begin{eqnarray*}
 A_n W_n^\Psi(t) &\equiv& A_n
 \int_0^t S_n(t-\tau)\Psi(\tau)dW(\tau) \\
 &=& J_n \int_0^t AS_n(t-\tau)\Psi(\tau)dW(\tau)
 = J_n \left[\int_0^t S_n(t-\tau)A\Psi(\tau)dW(\tau)\right].
\end{eqnarray*}
Analogously,
$$ A_n W^\Psi(t) =
  J_n \left[\int_0^t S(t-\tau)A\Psi(\tau)dW(\tau)\right].
$$
By (\ref{eq30}) we have
\begin{eqnarray*}
N_{n,1}(t) &=& |J_n\int_0^t [S_n(t-\tau)-S(t-\tau)]
  A\Psi(\tau)dW(\tau)|_H \\
  &\le & |\int_0^t [S_n(t-\tau)-S(t-\tau)]
  A\Psi(\tau)dW(\tau)|_H \;.
\end{eqnarray*}
From assumptions, $A\Psi \in \mathcal{N}^2(0,T;L_2^0)$. Then the
term $[S_n(t-\tau)-S(t-\tau)]A\Psi(\tau)$ may be treated like
the difference $\Phi_n$ defined by (\ref{eq21b}).

Hence, from (\ref{eq30}) and (\ref{eq27}), for the first term of
the right hand side of (\ref{eq28}) we have
$$ \lim_{n\to +\infty}\;\; \sup_{t\in [0,T]}
 \mathbb{E}(N_{n,1}^2 (t)) \to 0. $$
For the second term of (\ref{eq28}), that is $N_{n,2}^2(t)$, we can 
follow the same steps as above for proving (\ref{eq27}).
\begin{eqnarray*}
N_{n,2}(t) & = & |A_n W^\Psi (t) - A W^\Psi (t)|_H \\
 &\equiv & 
  \left| A_n \int_0^t S(t-\tau)\Psi (\tau)dW(\tau) 
       - A   \int_0^t S(t-\tau)\Psi (\tau)dW(\tau) \right|_H =  \\
 & = & \left|\int_0^t [A_n-A] S(t-\tau)\Psi (\tau)dW(\tau)\right|_H\;.
\end{eqnarray*}
 
From assumptions, $\Psi, A\Psi \in \mathcal{N}^2(0,T;L_2^0)$. 
Because $A_n, S(t), ~t\ge 0$ are bounded, then 
$A_n S(t-\cdot)\Psi(\cdot)
\in \mathcal{N}^2(0,T;L_2^0)$, too. Analogously, 
$AS(t-\cdot)\Psi(\cdot)=S(t-\cdot)A\Psi(\cdot)
\in \mathcal{N}^2(0,T;L_2^0)$.

Let us note that the set of all Hilbert-Schmidt operators acting from 
one separable Hilbert space into another one, equipped with the
operator norm defined on page \pageref{H-Sn} is a separable 
Hilbert space. Particularly, sum of two Hilbert-Schmidt operators
is a Hilbert-Schmidt operator, see e.g.\ \cite{Ba}. Therefore, 
we can deduce that the operator  $(A_n-A)$ $S(t-\cdot)\Psi(\cdot)\in 
\mathcal{N}^2(0,T;L_2^0)$, ~for any $~t\in [0,T]$. Hence, the term 
$[A_n-A]S(t-\tau)\Psi(\tau)$ may be treated like the difference 
$\Phi_n$ defined by (\ref{eq21b}). So, we obtain 
\begin{eqnarray*}
 \mathbb{E} \left( N_{n,2}^2(t)\right) &=& 
 \mathbb{E} \left( \int_0^t \left[ \sum_{j=1}^{+\infty} \left| 
 [A_n-A] S(t-\tau )\Psi (\tau)\,f_j \right|_H^2  \right]
 d\tau \right) \\
 &\le & \mathbb{E} \left(\sum_{j=1}^{+\infty} \int_0^T
 \left| [A_n-A] S(t-\tau)\Psi(\tau)\,f_j\right|_H^2 d\tau\right)
 < +\infty,
\end{eqnarray*}
for any $t\in [0,T]$.

By the convergence (\ref{eq31}), for any fixed $j$, 
$$ \int_0^T |[A_n-A]S(t-\tau) \Psi(\tau)\,f_j |_H^2 d\tau $$
tends to zero for $n\to +\infty$.

Summing up our considerations, we have 
$$ \lim_{n\to +\infty} \;\;\sup_{t\in [0,T]} \mathbb{E}
(N_{n,2}^2 (t)) \to 0\;.$$

So, we can deduce that
$$ \lim_{n\to +\infty} \;\;\sup_{t\in [0,T]} \mathbb{E}
  |A_n W_n^\Psi(t)-AW^\Psi(t)|_H^2=0 , $$
and then (\ref{deq16}) holds. \hfill $\blacksquare$

The main result of this section is the following.

\begin{theor} \label{coSW4}
Suppose that assumptions of Theorem \ref{pSW2} or Theorem
\ref{pSW2a} hold. Then the equation (\ref{eSW1}) has a strong
solution. Precisely, the convolution $W^\Psi$ defined by
(\ref{eSW18a}) is the strong solution to~(\ref{eSW1}) with
$X_0\equiv 0$.
\end{theor}
\bgproof In order to prove Theorem \ref{coSW4}, we have to show
only the condition (\ref{eSW3.1}). Let us note that the
convolution $W^\Psi (t)$ has integrable trajectories. Because the
closed unbounded linear 
operator $A$ becomes bounded on ($D(A),|\cdot|_{D(A)}$), see
\cite[Chapter 5]{We}, we obtain that $AW^\Psi (\cdot )\in
L^1([0,T];H)$, ~P-a.s. Hence, properties of convolution provide
integrability of the function $a(T-\tau)AW^\Psi (\tau)$ with
respect to $\tau$, what finishes the proof. \edproof

\end{document}